\pdfoutput=1    

\documentclass[11pt,a4paper]{article} 

\usepackage[OT2,T1]{fontenc}
\usepackage{amsmath}
\usepackage{amssymb,verbatim}
\usepackage{amsfonts}
\usepackage{amsthm}
\usepackage{pifont}
\usepackage{graphicx}
\usepackage{datetime}
\usepackage{scrdate}
\usepackage{enumitem}
\usepackage{geometry}
\usepackage{float}
\usepackage{bigints}
\usepackage{multirow}

\usepackage{hyperref}

\hypersetup{colorlinks=true, linkcolor=red, citecolor=red, pdfstartview=FitH, linktocpage=false}

\theoremstyle{definition}

\newcommand{\beq}{\begin{equation*}}
\newcommand{\eeq}{\end{equation*}}
\newcommand{\beqn}{\begin{equation}}
\newcommand{\eeqn}{\end{equation}}
\newcommand{\dd}{\mathrm{d}}
\newcommand{\ii}{\mathrm{i}}
\newcommand{\ee}{\mathrm{e}}

\newcommand{\ph}{\varphi}
\newcommand{\Ga}{\varGamma}

\begin{document}

\title{Simple evaluation of one of Malmst\'en's integrals}
\author{Uwe B\"asel}
\date{} 
\maketitle
\thispagestyle{empty}
\begin{abstract}
\noindent
The logarithmic integral no.\;4.325.7 from Gradshteyn and Ryzhik's tables of integrals was first evaluated by Malmst\'en. Recently, Blagouchine used contour integration methods to evaluate a family of logarithmic integrals that contains this integral. We evaluate the integral in a simple, straightforward manner mainly by means of real analysis. The main ingredients of the evaluation are the use of geometric series and Kummer's Fourier series expansion for the logarithm of the gamma function.
\\[0.2cm]
\\[0.2cm]
\end{abstract}

\section{Introduction}

In Gradshteyn and Ryzhik's tables of integrals \cite{Gradstein_Ryshik} one finds on p.\ 619 with no.\;4.325.7 the integral
\beqn \label{Eq:I}
  I(\ph)
= \int\limits_0^1 \frac{\ln\ln(1/x)\,\dd x}{1+2x\cos\ph+x^2}
= \int\limits_1^\infty \frac{\ln\ln x\,\dd x}{1+2x\cos\ph+x^2}
= \frac{\pi}{2\sin\ph}\ln\frac
	{(2\pi)^{\ph/\pi}\,\Ga\bigg(\dfrac{1}{2}+\dfrac{\ph}{2\pi}\bigg)}
	{\Ga\bigg(\dfrac{1}{2}-\dfrac{\ph}{2\pi}\bigg)}\,.
\eeqn
Blagouchine \cite[p.\ 76]{Blagouchine} (see also \cite{Malmsten_Wikipedia}) discovered that this integral was evaluated by Malmst\'en and published in \cite[p.\ 24]{Malmsten} already in the middle of the nineteenth century. In Blagouchine's paper it is Exercise 29 (h) on p.\ 76 to show that this integral can be deduced from a more general integral evaluation obtained by contour integration.

For $\ph=\pi/2$ we get
\beqn \label{Eq:Vardi-1b}
  I\!\left(\frac{\pi}{2}\right)
= \int_0^1 \frac{\ln\ln(1/x)\,\dd x}{1+x^2}
= \int_1^\infty \frac{\ln\ln u\,\dd u}{1+u^2}
= \frac{\pi}{2}\,\ln\frac{\Ga(3/4)\,\sqrt{2\pi}}{\Ga(1/4)}
\eeqn
which is no. 4.325.4 on p.\ 619 in \cite{Gradstein_Ryshik}. The substitution
\beq
  u = \tan y\,,\quad
  y = \arctan u\,,\quad
  \dd y = \frac{\dd u}{1+u^2}
\eeq
yields
\beqn \label{Eq:Vardi-1a}
  I\!\left(\frac{\pi}{2}\right)
= \int_{\pi/4}^{\pi/2} \ln\ln\tan y\,\dd y\,.  
\eeqn
This integral has no.\;4.229.7 on p.\ 577 in \cite{Gradstein_Ryshik}. It is called Vardi's integral \cite{Weisstein} since Vardi evaluated it in his paper \cite{Vardi} without knowing the author of this formula and the location of its proof. In his evaluation he used Dirichlet series, the Hurwitz $\zeta$-function and the Bohr-Mollerup Theorem. Vardi also considered the integrals
\beqn \label{Eq:Vardi2}
  I\!\left(\frac{\pi}{3}\right)
= \int_0^1 \frac{\ln\ln(1/x)\,\dd x}{1+x+x^2}
= \int_1^\infty \frac{\ln\ln x\,\dd x}{1+x+x^2}
= \frac{\pi}{\sqrt{3}}\,\ln\frac
	{\Ga(2/3)\,(2\pi)^{1/3}}{\Ga(1/3)}
\eeqn
and
\beqn \label{Eq:Vardi3}
\left.
\begin{aligned}
  I\!\left(\frac{2\pi}{3}\right)
= {} & \int_0^1 \frac{\ln\ln(1/x)\,\dd x}{1-x+x^2}
= \int_1^\infty \frac{\ln\ln x\,\dd x}{1-x+x^2}
= \frac{2\pi}{\sqrt{3}}\left[\frac{5}{6}\,\ln(2\pi)
- \ln\Ga\!\left(\frac{1}{6}\right)\right]\\
= {} & \frac{\pi}{\sqrt{3}}\,\ln\frac
	{\Ga(5/6)\,(2\pi)^{2/3}}{\Ga(1/6)}
\end{aligned}
\;\right\}
\eeqn
which are integrals no.\;4.325.5 and 4.325.6, respectively, in \cite[p.\ 619]{Gradstein_Ryshik}.
Integrals \eqref{Eq:Vardi-1b}, \eqref{Eq:Vardi-1a}, \eqref{Eq:Vardi2} and \eqref{Eq:Vardi3} were first evaluated by Malmst\'en and his colleagues (see \cite[pp.\ 22-23]{Blagouchine}).
(For further information on these and other logarithmic integrals see e.\,g. \cite{Adamchik}, \cite{Amdeberhan}, \cite{Blagouchine}, \cite{Boros_Moll}, \cite{Coeffey}, \cite{Malmsten}, \cite{Medina_Moll}, \cite{Moll}, \cite{Malmsten_Wikipedia}.)

In the following, we give a simple, straightforward proof for the evaluation of Malmst\'en's integral \eqref{Eq:I} mainly by means of real analysis without using other special functions than the gamma function.
The arguments used in this proof are self-evident (cp.\ the remarks concerning integral \eqref{Eq:Vardi-1a} in \cite[p.\ 308]{Vardi}, \cite[p.\ 47]{Amdeberhan} and \cite[pp.\ 1-2]{Moll}).

\section{Proof of Equation \eqref{Eq:I}}

We consider the denominator of the integrands in \eqref{Eq:I}. With $\cos\ph = (\ee^{\ii\ph}+\ee^{-\ii\ph})/2$ and applying Vieta's formula we have
\beq
  1 + 2x\cos\ph + x^2
= 1 + \left(\ee^{\ii\ph}+\ee^{-\ii\ph}\right)x + x^2
= \left(x + \ee^{\ii\ph}\right)\left(x + \ee^{-\ii\ph}\right)
= \left(1 + x\ee^{-\ii\ph}\right)\left(1 + x\ee^{\ii\ph}\right).
\eeq
With
\beq
  \left|x\ee^{-\ii\ph}\right|
= |x| \cdot \left|\ee^{-\ii\ph}\right|
= |x|
  \quad\mbox{and}\quad
  \left|x\ee^{\ii\ph}\right|
= |x| \cdot \left|\ee^{\ii\ph}\right|
= |x|
\eeq
we see that
\beq
  \left|x\ee^{-\ii\ph}\right|,\, \left|x\ee^{\ii\ph}\right| < 1
  \quad\mbox{if}\quad 0 \le x < 1\,.
\eeq
Expanding $\left(1 + x\ee^{-\ii\ph}\right)^{-1}$ and $\left(1 + x\ee^{\ii\ph}\right)^{-1}$ in geometric series yields the Cauchy product
\begin{align*}
  \frac{1}{\left(1 + x\ee^{-\ii\ph}\right)\left(1 + x\ee^{\ii\ph}\right)}
= \sum_{n=0}^\infty (-1)^n \left(x\ee^{-\ii\ph}\right)^n \cdot
  \sum_{n=0}^\infty (-1)^n \left(x\ee^{\ii\ph}\right)^n\\
& \hspace{-10.4cm} = \left(1 - x\ee^{-\ii\ph} + x^2\ee^{-2\ii\ph}
  - x^3\ee^{-3\ii\ph} +- \ldots\right)
  \left(1 - x\ee^{\ii\ph} + x^2\ee^{2\ii\ph}
  - x^3\ee^{3\ii\ph} +- \ldots\right)\\
& \hspace{-10.4cm} = 1 - \left(\ee^{\ii\ph}+\ee^{-\ii\ph}\right)x
  + \left(\ee^{2\ii\ph}+1+\ee^{-2\ii\ph}\right)x^2
  - \left(\ee^{3\ii\ph}+\ee^{\ii\ph}+\ee^{-\ii\ph}
	+\ee^{-3\ii\ph}\right)x^3
  +- \ldots\,.
\end{align*}
We denote by $a_n$ the coefficient of $x^n$ without the sign. So we have
\begin{align*}
  a_n
= {} & \sum_{k=0}^n \ee^{(n-2k)\ii\ph}
= \ee^{n\ii\ph}\,\sum_{k=0}^n \ee^{-2k\ii\ph}
= \ee^{n\ii\ph}\,\sum_{k=0}^n \left(\ee^{-2\ii\ph}\right)^k
= \ee^{n\ii\ph}\,\frac{1-\ee^{-2(n+1)\ii\ph}}{1-\ee^{-2\ii\ph}}\\
= {} & \frac{\ee^{(n+1)\ii\ph}}{\ee^{\ii\ph}}\,
  \frac{1-\ee^{-2(n+1)\ii\ph}}{1-\ee^{-2\ii\ph}}
= \frac{\ee^{(n+1)\ii\ph}-\ee^{-(n+1)\ii\ph}}{\ee^{\ii\ph}-\ee^{-\ii\ph}}
= \frac{\sin(n+1)\ph}{\sin\ph}\,.
\end{align*}
So the first integral in \eqref{Eq:I} may be written as
\begin{align*}
  I(\ph)
= {} & \int_0^1 \ln\left(\ln\frac{1}{x}\right)\,
  \sum_{n=0}^\infty (-1)^n\,\frac{\sin(n+1)\ph}{\sin\ph}\,x^n\,\dd x\displaybreak[0]\\
= {} & \frac{1}{\sin\ph}\sum_{n=0}^\infty (-1)^n\sin(n+1)\ph
  \underbrace{\int_0^1 x^n\ln\left(\ln\frac{1}{x}\right)\dd x}_{J_n}\,.
\end{align*}
The substitution $u=\ln(1/x)$ gives
\beq
  J_n
= \int_0^\infty \ee^{-(n+1)u}\ln u\,\dd u\,.
\eeq
Now we substitute $y = (n+1)u$ and get
\begin{align*}
  J_n
= {} & \frac{1}{n+1}\int_0^\infty \ee^{-y}\ln\frac{y}{n+1}\:\dd y
= \frac{1}{n+1}\left[\int_0^\infty \ee^{-y}\ln y\,\dd y
- \ln(n+1)\int_0^\infty \ee^{-y}\,\dd y\right]\\[0.1cm]
= {} & \frac{1}{n+1}\left[\int_0^\infty \ee^{-y}\ln y\,\dd y
- \ln(n+1)\right].
\end{align*}
We have
\beq
  \int_0^\infty \ee^{-y}\ln y\,\dd y = -\gamma\,,
\eeq 
where $\gamma = 0.5772156649\ldots$ is the Euler-Mascheroni constant (see e.\,g.\ \cite[Vol.\ II, pp.\ 801, 813]{Fichtenholz}, \cite[pp.~176-177]{Boros_Moll}), hence
\beq
  J_n
= -\frac{\gamma + \ln(n+1)}{n+1}\,.
\eeq
It follows that
\begin{align*}
  I(\ph)
= {} & -\!\frac{1}{\sin\ph}\,\sum_{n=0}^\infty (-1)^n\,\frac{\gamma
+ \ln(n+1)}{n+1}\,\sin(n+1)\ph\\
= {} & -\!\frac{1}{\sin\ph}\left[\gamma\,\sum_{n=0}^\infty
  (-1)^n\,\frac{\sin(n+1)\ph}{n+1} 
+ \sum_{n=0}^\infty (-1)^n\,\frac{\ln(n+1)\,\sin(n+1)\ph}{n+1}
  \right]\\
= {} & \frac{1}{\sin\ph}\left[-\gamma\,\sum_{n=0}^\infty
  (-1)^n\,\frac{\sin(n+1)\ph}{n+1} 
+ \sum_{n=0}^\infty (-1)^{n+1}\,\frac{\ln(n+1)\,\sin(n+1)\ph}{n+1}
  \right]\\
= {} & \frac{1}{\sin\ph}\left[-\gamma\,\sum_{n=1}^\infty
  (-1)^{n+1}\,\frac{\sin n\ph}{n} 
+ \sum_{n=2}^\infty (-1)^n\,\frac{\ln n\,\sin n\ph}{n}\right].
\end{align*}
With 
\beq
  \sum_{n=1}^\infty (-1)^{n+1}\,\frac{\sin n\ph}{n}
= \frac{\ph}{2}\,,\quad -\pi < \ph < \pi\,, 
  \quad\mbox{\cite[p.\,387]{Knopp}}\quad
\eeq
we have
\beqn \label{Eq:I2}
  I(\ph)
= \frac{1}{\sin\ph}\left[-\frac{\gamma\ph}{2} 
+ \sum_{n=2}^\infty (-1)^n\,\frac{\ln n\,\sin n\ph}{n}\right]
\eeqn
so that it remains to evaluate the series in \eqref{Eq:I2}. To this end, we apply Kummer's Fourier series expansion for the logarithm of the gamma function (see e.\,g.\ \cite[Vol.\ III, pp.\ 452-455]{Fichtenholz} with Kummer's proof). An equivalent but simplified version (as one may prove) is provided in \cite{Gamma_Wikipedia}: 
\beqn \label{Eq:Gamma}
  \ln\Ga(x)
= \left(\frac{1}{2}-x\right)(\gamma+\ln 2)+(1-x)\ln\pi-\frac{1}{2}\ln\sin\pi x
+ \frac{1}{\pi}\sum_{n=1}^{\infty}\frac{\ln n \cdot \sin 2\pi nx}{n}\,,\; 0<x<1\,.
\eeqn
One sees that the series in \eqref{Eq:I2} is alternating while the series in \eqref{Eq:Gamma} is not. In order to remove this problem we observe that
\beq
  \sin(n\pi-n\ph)
= \sin n\pi \cos n\ph - \cos n\pi\sin n\ph
= -\cos n\pi \sin n\ph
\eeq
with
\beq
  \cos n\pi
= \left\{\begin{array}{rl}
	 1 & \mbox{if\;\: $n$ is even}\,,\\[0.1cm]
	-1 & \mbox{if\;\: $n$ is odd}\,.  
  \end{array}\right.
\eeq
So we have
\beq
  (-1)^n \sin n\ph = -\sin n(\pi-\ph)\,,
\eeq
and \eqref{Eq:I2} becomes
\beqn \label{Eq:I3}
  I(\ph)
= -\frac{1}{\sin\ph}\left[\frac{\gamma\ph}{2} 
+ \sum_{n=2}^\infty \frac{\ln n\,\sin n(\pi-\ph)}{n}\right].
\eeqn
The substitution 
\beq
  2\pi x = \pi - \ph \quad\Longleftrightarrow\quad
  x = \frac{1}{2} - \frac{\ph}{2\pi}
\eeq
in \eqref{Eq:Gamma} leads to
\beq
  \sum_{n=2}^\infty \frac{\ln n\,\sin n(\pi-\ph)}{n}
= \pi\ln\Ga\left(\frac{1}{2}-\frac{\ph}{2\pi}\right)
- \frac{\ph}{2}\,(\gamma+\ln 2\pi) - \frac{\pi}{2}\,\ln\pi
+ \frac{\pi}{2}\,\ln\cos\frac{\ph}{2} 
\eeq
for $-\pi < \ph < \pi$ (cp.\ \cite[p.\ 67, Exercise 20 (b)]{Blagouchine} and \cite{Malmsten_Wikipedia}).
It follows that
\begin{align*}
  I(\ph)
= {} & -\!\frac{1}{\sin\ph}
  \left[\pi\ln\Ga\!\left(\frac{1}{2}-\frac{\ph}{2\pi}\right)
- \frac{\ph}{2}\,\ln(2\pi) 
- \frac{\pi}{2}\ln\pi
+ \frac{\pi}{2}\,\ln\cos\frac{\ph}{2}
\right]\\
= {} & \frac{\pi}{\sin\ph}
  \left[-\!\ln\Ga\!\left(\frac{1}{2}-\frac{\ph}{2\pi}\right)
+ \frac{\ph}{2\pi}\,\ln(2\pi) 
+ \frac{1}{2}\ln\pi
- \frac{1}{2}\,\ln\cos\frac{\ph}{2}
\right]\\
= {} & \frac{\pi}{\sin\ph}\,
  \ln\frac
	{(2\pi)^{\frac{\ph}{2\pi}}\,\sqrt{\pi}}
	{\Ga\!\left(\frac{1}{2}-\frac{\ph}{2\pi}\right)
		\sqrt{\cos(\ph/2)}}\,.
\end{align*}
Applying Euler's reflection formula gives
\beq
  \Ga\!\left(\frac{1}{2}-\frac{\ph}{2\pi}\right)
  \Ga\!\left(\frac{1}{2}+\frac{\ph}{2\pi}\right)
= \frac
	{\pi}
	{\sin\left(\pi\left(\frac{1}{2}-\frac{\ph}{2\pi}\right)\right)}
= \frac{\pi}{\cos(\ph/2)}\,,
\eeq
and so we may write our result as
\beq
  I(\ph)
= \frac{\pi}{2\sin\ph}\,
  \ln\frac
	{(2\pi)^{\ph/\pi}\,\Ga\bigg(\dfrac{1}{2}+\dfrac{\ph}{2\pi}\bigg)}
	{\Ga\bigg(\dfrac{1}{2}-\dfrac{\ph}{2\pi}\bigg)}\,,\quad
  -\pi < \ph < \pi\,.  
\eeq

\bigskip\noindent
{\bf Uwe B\"asel}, Hochschule f\"ur Technik, Wirtschaft und Kultur Leipzig, Fakult\"at f\"ur Maschinenbau und Energietechnik, PF 30\,11\,66, 04251 Leipzig, Germany, \texttt{uwe.baesel@htwk-leipzig.de}

\end{document}